\documentclass{amsart}[11pt]
\usepackage{amsmath, amsthm, amsfonts}

\newcommand{\dd}{\partial}
\newcommand{\SSS}{{\cal S}}
\newcommand{\sss}{{\frak s}}
\newcommand{\cal}{\mathcal}

\swapnumbers
\newtheorem{theorem}{Theorem}[section]
\newtheorem{lemma}[theorem]{Lemma}
\newtheorem{corollary}[theorem]{Corollary}
\newtheorem{definition}[theorem]{Definition}
\newtheorem{proposition}[theorem]{Proposition}
\newtheorem*{maintheorem}{Main Theorem}
\newtheorem{observation}[theorem]{Observation}
\newtheorem*{remark}{Remark}

\newcommand{\R}{\mathbb R}
\newcommand{\Z}{\mathbb Z}

\newcommand{\iso}{\cong}
\newcommand{\im}{\text{im\ }}
\newcommand{\comment}[1]{}
\renewcommand{\tilde}[1]{\widetilde{#1}}

\newcommand{\ssection}[1]{\section{#1}\setcounter{equation}{0}}

 \def\Slash#1{
  \begin{picture}(5,6)(0,0)
  \put(-.7,-1.2){\line(5,6)5}
  \end{picture}
  \kern-.8em#1}
 \def\slash#1{
  \begin{picture}(5,6)(0,0)
  \put(-1.5,-1.7){\line(5,6)5}
  \end{picture}
  \kern-.8em#1}

\def\sd{\Slash \partial}

\begin{document}

\title{Conformally invariant powers of the ambient Dirac operator}
\author[J. Holland]{Jonathan Earl Holland}
\address{Department of Mathematics, Princeton University, Princeton, New Jersey 08544}
\email{jholland@math.princeton.edu}

\author[G. Sparling]{George A. J. Sparling}
\address{Department of Mathematics, University of Pittsburgh, Pittsburgh, Pennsylvania 15260}
\email{sparling@math.pitt.edu}




\begin{abstract}
This paper constructs a family of conformally invariant differential operators acting on spinor densities with leading part a power of the Dirac operator.  The construction applies for all powers in odd dimensions, and only for finitely many powers in even dimensions.  These operators arise naturally as obstructions to formal solution of the Dirac equation on the Fefferman-Graham ambient space with prescribed boundary conditions.
\end{abstract}

\maketitle

\ssection{Introduction}\label{intro}

Recent work in even-dimensional conformal geometry \cite{Branson},\cite{FeffGra2}, \cite{Gov}, \cite{GraZ} has revealed the importance of conformally invariant powers of the Laplacian on a conformal manifold; that is, of operators $P_k$ whose principal part is the same as $\Delta^k$ with respect to a representative of the conformal structure.  These invariant powers of the Laplacian were first defined in \cite{GJMS} in terms of the Fefferman-Graham \cite{FeffGra1} ambient Lorentzian structure.  The present paper defines conformally invariant powers of the Dirac operator by generalizing the construction of \cite{GJMS}.

It is well-known that the usual Dirac operator $\sd$, under a conformal transformation $g\mapsto \widehat{g}=\lambda^2 g$, transforms (acting on spinor densities of weight $0$) via
$$\widehat{\sd}\psi={1\over\lambda^2}\sd(\lambda\psi).$$
Such an operator is called conformally covariant of bidegree $(-2,1)$ in \cite{Branson}.  However, the operator we shall study in this paper is not the usual Dirac operator, but rather comes from the Dirac operator of the ambient Lorentzian manifold constructed in \cite{FeffGra1}.  Our approach is nearly parallel to the construction of conformally invariant powers of the Laplacian in \cite{GJMS}.  However, there are many subtleties involved in the study of the Dirac operator which are not present in the case of the Laplacian.

Our main results are summarized as follows, stated somewhat imprecisely:
\begin{maintheorem} Let $M$ be a manifold of dimension $n$.  Let $G, \tilde{G}$ be as in {\rm \cite{GJMS}}.  
\begin{enumerate}
\item Let $\Psi$ be a spinor on $G$ of conformal weight $w=p-{1\over 2}n$. Let $\psi$ be an extension of $\Psi$ to a spinor on $\tilde{G}$ also of conformal weight $w$.  Then $\bigl(x{\sd}^{2p}\psi\bigr)|_G$ is a spinor on $G$ of conformal weight $-p-{1\over 2}n+1$ independent of the choice of extension $\psi$.  Thus $\Psi\mapsto\bigl(x{\sd}^{2p}\psi\bigr)|_G$ is a conformally invariant operator.
\item Let $\Psi$ have conformal weight $w=p-{1\over 2}n+1$.  Let $\psi$ be a {\rm ``preferred''} extension of $\Psi$ of conformal weight $w$.  Then $\bigl({\sd}^{2p+1}\psi\bigr)|_G$ is a spinor on $G$ of conformal weight $-p-{1\over 2}n$ independent of the choice of preferred extension $\psi$.  Thus $\Psi\mapsto\bigl({\sd}^{2p+1}\psi\bigr)|_G$ is a conformally invariant operator.
\end{enumerate}
\end{maintheorem}
Here $x$ is a $C^\infty$-linear operator related to the Euler field, which will be defined in Section \ref{gradedLie}.  Part 1 of the main theorem is not surprising, considering the result Proposition 2.1 in \cite{GJMS}, and we shall show that $\bigl(x{\sd}^{2p}\psi\bigr)|_G$ arises naturally as the obstruction to an asymptotic expansion problem by analogy with Proposition 2.2 in \cite{GJMS}.  Part 2 has no analog in \cite{GJMS}.  Furthermore, we show that the operator $\bigl({\sd}^{2p+1}\psi\bigr)|_G$ arises naturally as the obstruction to a different asymptotic expansion problem.

\ssection{The ambient space}\label{ambientspace}
Let $M$ be a manifold of dimension $n$ with a conformal metric $[g]$ of signature $(p,q)$.  Let ${\cal E}[w]$ be the line bundle of conformal densities of weight $w$ on $M$.  Let $G\xrightarrow{\pi} M$ be the principal $\R_{>0}$-bundle of representatives of the conformal class on $M$.  Define a one-parameter family of dilation diffeomorphisms $\delta_s:G\rightarrow G$ by $\delta_s(x,g)=(x,s^2g)$.  The {\it Euler field} is then the vector field $X$ on $G$ defined by
$$Xf(p)={d\over ds}f(\delta_sp)|_{s=0},$$
for $f\in C^\infty(G), p\in G$.  In abstract indices, we shall denote the Euler field by $x^a$.  There is a tautological bilinear form $g_0$ on $G$ of signature $(p,q,0)$ given at a point $(x,g)\in G$ by
$$g_0(Y,Z)=g(\pi_*Y,\pi_*Z).$$
Note that
\begin{equation}\label{x|_G} g_0(X,-)=0\end{equation}

Let $\tilde{G}$ be a formal neighborhood of $G\times\{0\}$ in $G\times [0,1)$.  The $\delta_s$ extends to a diffeomorphism of $\tilde{G}$, and so $X$ extends to a vector field on $\tilde{G}$ as well.  In \cite{FeffGra1} it is proven that $g_0$ has a unique formal extension to a formal metric $\tilde{g}$ of signature $(p+1,q+1)$ on $\tilde{G}$ such that $\delta_s^*\tilde{g}=s^2\tilde{g}$ and $\hbox{Ric}(\tilde{g})=0$.  For odd $n$, this formal extension is defined to infinite order.  For even $n$, the extension is defined in general only to order $n/2$.

We suppose that $M$ is a spin manifold (i.e., $w_2(M)=0$).  Then since $G$ and $G\times [0,1)$ are retractable to $M$, $G$ and $\tilde{G}$ each carry a spin structure as well.  Let $\SSS(M),\SSS(G),\SSS(\tilde{G})$ be the spin bundles of $M$, $G$ and $\tilde{G}$, respectively.  Let $\SSS(G)[w]$ (resp. $\SSS(\tilde{G})[w]$) denote the sheaf of all sections $\psi$ of $\SSS(G)$ (resp. $\SSS(\tilde{G})|_G$) such that ${\cal L}_X \psi=w\psi$.  

\begin{observation} \label{SGstruct1} There is a natural inclusion $\SSS(G)\subset\SSS(\tilde{G})|_G$.  Consider the relationship between $\SSS(M)$ and $\SSS(G)$.  Since the metric on $G$ is degenerate of signature $(p,q,0)$, its Clifford algebra is isomorphic to $Cl(p,q)$ with one element $x$ adjoined which anticommutes with every other element and such that $x^2=0$.  Thus, in any irreducible representation of $Cl(p,q,0)$, one has $x=0$.  Hence in particular 
\begin{equation}\label{SGstruct2}
\SSS(G)\iso\SSS(M)\times_M G.
\end{equation}
\end{observation}

Henceforth we shall work primarily with the ambient structure on $\tilde{G}$.  Thus denote by $g_{ab}$ the Fefferman-Graham metric on $\tilde{G}$.  Let $\dd_a$ be the Levi-Civita connection for the ambient structure.  Let us record the following facts from \cite{FeffGra1} and \cite{GJMS} for future reference
\begin{equation}\label{gab}
\dd_ax_b=g_{ab}
\end{equation}
and
\begin{equation}\label{curv}
\dd_b\dd_cx_d=R_{abcd}x^a=0.
\end{equation}

\ssection{The graded Lie algebra}\label{gradedLie}

Let $\gamma:T_p\tilde{G}\rightarrow\hbox{End}\bigl(\SSS(\tilde{G})_p\bigr)$ be the Clifford representation of the tangent space at $p$.  We may regard $\gamma$ as an element of $T^*_p\tilde{G}\otimes\SSS(\tilde{G})$, and thus identify $\gamma$ with the indexed operator $\gamma_a$.  The Clifford algebra $Cl(\tilde{G})$ is defined by the relation
\begin{equation}\label{cl}
\gamma_a\gamma_b+\gamma_b\gamma_a=2g_{ab}
\end{equation}
and give this algebra the usual $\Z_2$-grading by the number of factors of $\gamma$ occuring in its elements.  Let $x^a$ denote the Euler field.  Consider the graded Lie algebra $\frak g$ generated by 
\begin{equation}\begin{aligned}
x&=\gamma^ax_a\\
y&=\gamma^a\dd_a\\
h&=x^a\dd_a+{n+2\over 2},\end{aligned}\end{equation}
where the grading is again determined by parity of the number of $\gamma$'s in an element.  Recall that in a $\Z_2$-graded Lie algebra, the commutator of two elements $\alpha$ of degree $p$ and $\beta$ of degree $q$ is given in the enveloping algebra $E({\frak g})$ by
$$[\alpha,\beta]=\alpha\beta-(-1)^{pq}\beta\alpha.$$
Let $Q=x^2$, $\Delta=y^2$ in $E({\frak g})$.  Note that $y$ is precisely the Dirac operator $\sd$.  Of course, by the Lichnerowicz formula, $y^2={\sd}^2=\Delta+{1\over 4}\hbox{tr\ Ric}(\tilde{g})=\Delta$ (in our case) where $\Delta$ is the usual covariant Laplacian acting on spinors.  Furthermore, it is clear that $Q=0$ on $G$.

\begin{proposition}  The following commutation relations hold in $\frak g$:
\begin{eqnarray}
\label{xx}
[x,x]&=&2Q\\
\label{yy}
[y,y]&=&2\Delta\\
\label{xy}
[x,y]&=&2h\\
\label{Qx}
[Q,x]&=&0\\
\label{Qy}
[Q,y]&=&-2x\\
\label{Qh}
[Q,h]&=&-2Q\\
\label{Dx}
[\Delta,x]&=&2y\\
\label{Dy}
[\Delta,y]&=&0\\
\label{Dh}
[\Delta,h]&=&2\Delta\\
\label{DQ}
[\Delta,Q]&=&4h\\
\label{xh}
[x,h]&=&-x\\
\label{yh}
[y,h]&=&y
\end{eqnarray}
\end{proposition}

\begin{remark}  $\frak g$ is isomorphic to the orthosymplectic Lie superalgebra ${\frak o}{\frak s}{\frak p}_{1,2}$, although this fact is not used anywhere in the paper.
\end{remark}

\begin{proof}[Proof of Propostion.]  \eqref{xx}, \eqref{yy} are trivial consequences of the definition of $Q$ and $\Delta$, respectively.  

\eqref{xy}: For a spinor $\psi$, we have
\begin{align*}
[x,y]\psi&=\gamma^ax_a\gamma^b\dd_b\psi+\gamma^b\dd_b\gamma^ax_a\psi &&\\
&=\gamma^a\gamma^bx_a\dd_b\psi+\gamma^b\gamma^a(\dd_bx_a)\psi+\gamma^b\gamma^ax_a\dd_b\psi &&\\
&=(\gamma^a\gamma^b+\gamma^b\gamma^a)x_a\dd_b\psi+\gamma^b\gamma^a(\dd_bx_a)\psi &&\\
&=2g^{ab}x_a\dd_b\psi+\gamma^b\gamma^ag_{ba}\psi &&\text{by \eqref{cl} and \eqref{gab}}\\
&=2h\psi
\end{align*}

\eqref{Qx}:  Follows immediately from $Q=x^2$.

\eqref{Qy}:  
\begin{align*}
[Q,y]\psi&=Q\gamma^a\dd_a\psi-\gamma^a\dd_a(Q\psi) &&\\
&=Q\gamma^a\dd_a\psi-\gamma^aQ\dd_a\psi-\gamma^a(\dd_aQ)\psi&&\\
&=-\gamma^a\dd_a(x^bx_b)\psi&&\text{by definition of $Q$}\\
&=-2\gamma^ax^b(\dd_ax_b)\psi&&\\
&=-2\gamma^ax^bg_{ab}\psi&&\text{by \eqref{gab}}\\
&=-2x\psi
\end{align*}

\eqref{Qh}:  
\begin{align*}
[Q,h]\psi&=Qx^a\dd_a\psi-x^aQ\dd_a\psi-x^a(\dd_aQ)\psi&&\\
&=-2x^ax^b(\dd_ax_b)\psi&&\\
&=-2Q\psi&&\text{by \eqref{gab}}
\end{align*}

\eqref{Dx}:  
\begin{align*}
[\Delta,x]\psi&=\Delta(x\psi)-x(\Delta\psi)&&\\
&=(\Delta x)\psi+2(\dd^ax)(\dd_a\psi)+x(\Delta\psi)-x(\Delta\psi)&&\\
&=\gamma^a(\dd^b\dd_bx_a)\psi+2\gamma^a\dd_a\psi &&\text{by \eqref{gab}}\\
&=\gamma^a(\dd^bg_{ba})\psi+2y\psi &&\text{by \eqref{gab}}\\
&=2y\psi &&\text{since $\dd^b$ is Levi-Civita}
\end{align*}

\eqref{Dy}:  Follows from $\Delta=y^2$.

\eqref{Dh}:  For this, note that the curvature endomorphism $[\dd_a,\dd_b]$ extends to an endomorphism of the spin bundle, which we denote by $R_{ab}$.  It then follows from \eqref{curv} that
\begin{equation}\label{curv2}
x^aR_{ab}=0.
\end{equation}
Differentiating this relation gives
\begin{align*}
0&=\dd^bx^aR_{ab}\psi&&\\
&=(\dd^bx^a)R_{ab}\psi+x^a\dd^bR_{ab}\psi&&\\
&=x^a\dd^bR_{ab}\psi&&\text{by skew-symmetry of $R_{ab}$.}
\end{align*}
Thus we have
\begin{equation}\label{curv3}
x^a\dd^bR_{ab}=0
\end{equation}
Applying \eqref{curv2} and \eqref{curv3} respectively yields the following identity:
\begin{equation}\label{curv4}
x^a\dd_a\dd^b\dd_b\psi=x^a\dd^b\dd_b\dd_a\psi
\end{equation}

Now, we have
\begin{align*}
[\Delta,h]\psi&=\dd_b\dd^bx^a\dd_a\psi-x^a\dd_a\dd_b\dd^b\psi&&\\
&=\dd_b(\dd^bx^a)\dd_a\psi+(\dd^bx^a)\dd_b\dd_a\psi+x^a\dd^b\dd_b\dd_a\psi-x^a\dd_a\dd_b\dd^b\psi&&\\
&=\dd_bg^{ba}\dd_a\psi + g^{ba}\dd_b\dd_a\psi+x^a\dd^b\dd_b\dd_a\psi-x^a\dd_a\dd^b\dd_b\psi&&\text{by \eqref{gab}}\\
&=2\Delta\psi+x^a\dd^b\dd_b\dd_a\psi-x^a\dd_a\dd^b\dd_b\psi&&\\
&=2\Delta\psi&&\text{by \eqref{curv4}}\\
\end{align*}

\eqref{DQ}: 
\begin{align*}
[\Delta,Q]\psi&=(\Delta Q)\psi+2(\dd_a Q)(\dd^a\psi)&&\\
&=2(\dd^bx^a)(\dd_bx_a)\psi+x^a(\dd^b\dd_bx_a)\psi+4(x^c\dd_ax_c)(\dd^a\psi)&&\\&=(2(n+2)+4x_a\dd^a)\psi&&\text{by \eqref{gab}}\\
&=4h&&
\end{align*}

\eqref{xh}:  
\begin{align*}
[x,h]\psi&=\gamma^ax_ax^b\dd_b\psi-x^b\dd_b\gamma^ax_a&&\\
&=\gamma^ax_ax^b\dd_b\psi-x^b\gamma^a(\dd_bx_a)\psi-x^a\gamma^ax_a\dd_b\psi&&\\
&=-x^b\gamma^a(\dd_bx_a)&&\\
&=-x&&\text{by \eqref{gab}}
\end{align*}

\eqref{yh}:  
\begin{align*}
[y,h]&=\gamma^a\dd_ax^b\dd_b-x^b\dd_b\gamma^a\dd_a&&\\
&=\gamma^a(\dd_ax^b)\dd_b+\gamma^ax^b\dd_a\dd_b-x^b\gamma^a\dd_b\dd_a&&\\
&=\gamma^b\dd_b +\gamma^ax^bR_{ab}&&\text{by \eqref{gab}}\\
&=\gamma^b\dd_b&&\text{by \eqref{curv3}}\\
&=y&&
\end{align*}
\end{proof}

The opposite algebra ${\frak g}^\circ$ is defined to be the graded Lie algebra generated by the elements of $\frak g$ with the relation 
$$[\alpha,\beta]_{{\frak g}^\circ}=(-1)^{pq}[\alpha,\beta]_{\frak g},$$
where $\alpha$ and $\beta$ are homogeneous of $\Z_2$-degree $p$ and $q$, respectively.  We have the following corollary to Proposition 1, easily verified by inspection of the relations for $\frak g$:

\begin{corollary}\label{interchange}  The operation interchanging $x$ with $-y$, $Q$ with $-\Delta$, and $h$ with $-h$ determines an isomorphism
$${\frak g}\iso {\frak g}^\circ.$$
\end{corollary}

\begin{proposition}  Let $q(t)$ be a polynomial in one indeterminant and $p$ a nonnegative integer.  Then the following identities hold in $E({\frak g})$:
\begin{eqnarray}
\label{qhyp}
q(h)y^p&=y^pq(h-p)\\
\label{ypqh}
y^pq(h)&=q(h+p)y^p\\
\label{qhxp}
q(h)x^p&=x^pq(h+p)\\
\label{xpqh}
x^pq(h)&=q(h-p)x^p
\end{eqnarray}
\end{proposition}
\begin{proof}  It suffices to show \eqref{qhyp}, since \eqref{ypqh} is then immediate, and \eqref{qhxp}, \eqref{xpqh} follow by Corollary \ref{interchange}.  By linearity, it suffices to check the case where $q(h)=h^r$.  Furthermore, it suffices to prove the result for $r=1$ since then
$$h^{r-1}hy^p=h^{r-1}y^p(h-p),$$ 
and the result follows by induction.  Now,
\begin{align*}
hy^p&=[h,y]y^{p-1}+yhy^{p-1}&&\\
&=-y^p+yhy^{p-1}&&\text{by \eqref{yh}}
\end{align*}
and the result now follows by induction on $p$.
\end{proof}

\begin{proposition}\label{commutators} Let $p$ be a positive integer.  The following equations hold in $E({\frak g})$:
\begin{eqnarray}
\label{y2px}
[y^{2p},x]&=&2py^{2p-1}\\
\label{y2p+1x}
[y^{2p+1},x]&=&2y^{2p}(h-p)\\
\label{y2px2}
[y^{2p},x^2]&=&4py^{2p-2}(h-p+1)\\
\label{y2p+1x2}
[y^{2p+1},x^2]&=&4py^{2p-1}(h-p)+2y^{2p}x\\
\label{x2py}
[x^{2p},y]&=&2px^{2p-1}\\
\label{x2p+1y}
[x^{2p+1},y]&=&2x^{2p}(h+p)\\
\label{x2py2}
[x^{2p},y^2]&=&-4px^{2p-2}(h+p-1)\\
\label{x2p+1y2}
[x^{2p+1},y^2]&=&2x^{2p}y-4px^{2p-1}(h+p)
\end{eqnarray}
\end{proposition}

\begin{proof}  Note that it suffices to prove \eqref{y2px}, \eqref{y2p+1x}, \eqref{y2px2}, and \eqref{y2p+1x2}, for then the remaining identities follow by Corollary \ref{interchange}.

We first show \eqref{y2px}.  Write
\begin{align*}
[y^{2p},x]&=\sum_{j=1}^{2p}(-1)^{j-1}y^{2p-j}[x,y]y^{j-1}&&\\
&=2\sum_{j=1}^{2p} (-1)^{j-1}y^{2p-j}hy^{j-1}&&\\
&=2\sum_{j=1}^{2p} (-1)^{j-1}y^{2p-1}(h-j+1)&&\text{by \eqref{qhyp}}\\
&=2y^{2p-1}\sum_{j=1}^{2p}(-1)^j(j-1)&&\\
&=2py^{2p-1}
\end{align*}

Now, for \eqref{y2p+1x},
\begin{align*}
[y^{2p+1},x]&=\sum_{j=1}^{2p+1}(-1)^{j-1}y^{2p+1-j}[x,y]y^{j-1}&&\\
&=2\sum_{j=1}^{2p+1}(-1)^{j-1}y^{2p}(h-j+1)&&\text{by \eqref{qhyp}}\\
&=2y^{2p}\bigl(h-\sum_{j=1}^{2p+1}(-1)^{j-1}(j-1)\bigr)\\
&=2y^{2p}(h-p)
\end{align*}

\eqref{y2px2} is proven as follows
\begin{align*}
[y^{2p},x^2]&=[y^{2p},x]x+x[y^{2p},x]&&\\
&=2p(y^{2p-1}x+xy^{2p-1})&&\text{by \eqref{y2px}}\\
&=2p[y^{2p-1},x]&&\\
&=4py^{2p-2}(h-p+1)&&\text{by \eqref{y2p+1x}}
\end{align*}

Finally, for \eqref{y2p+1x2},
\begin{align*}
[y^{2p+1},x^2]&=[y^{2p+1},x]x-x[y^{2p+1},x]&&\\
&=2y^{2p}(h-p)x-2xy^{2p}(h-p)&&\text{by \eqref{y2p+1x}}\\
&=2y^{2p}x(h-p+1)-2xy^{2p}(h-p)&&\text{by \eqref{qhxp}}\\
&=2[y^{2p},x](h-p)+2y^{2p}x&&\\
&=4py^{2p-1}(h-p)+2y^{2p}x&&\text{by \eqref{y2px}}
\end{align*}
\end{proof}

\begin{proposition} Let $p$ be a positive integer.  For an indeterminant $t$, let $[t]^p$ denote the following polynomial
$$[t]^p=t(t+1)\dots (t+p-1).$$
Then the following equations hold in $E({\frak g})$:
\begin{eqnarray}
\label{y2p}
y^{2p}x^{2p}&=&2^{2p}p![h]^p+x^2Z_{2p}\\
\label{y2p+1}
y^{2p+1}x^{2p+1}&=&-2^{2p+1}p![h-1]^{p+1}+ 2^{2p}p![h]^pyx+x^2Z_{2p+1}
\end{eqnarray}
for some polynomials $Z_{2p}, Z_{2p+1}$ in $x,y,h$.
\end{proposition}

\begin{proof}  We first prove \eqref{y2p}.  We have
\begin{align*}
y^{2p}x^{2p}&=[y^{2p},x^2]x^{2p-2}+x^2y^{2p}x^{2p-2}&&\\
&=4py^{2p-2}(h-p+1)x^{2p-2}+O(x^2)&&\text{by \eqref{y2px2}}\\
&=4py^{2p-2}x^{2p-2}(h+p-1)+O(x^2)&&\text{by \eqref{qhxp}}
\end{align*}
and \eqref{y2p} follows by induction.

For \eqref{y2p+1},
\begin{align*}
y^{2p+1}x^{2p+1}&=-xy^{2p+1}x^{2p}+[y^{2p+1},x]x^{2p}&&\\
&=-xy^{2p+1}x^{2p}+2x^2y^{2p}x^{2p-2}(h+p)&&\\
&\quad +2[y^{2p},x^2]x^{2p-2}(h+p)&&\text{by \eqref{y2p+1x}, \eqref{qhxp}}\\
&=-xy^{2p+1}x^{2p}+8py^{2p-2}x^{2p-2}(h+p-1)+O(x^2)&&\text{by \eqref{y2px2}, \eqref{qhxp}}\\
&=-xy2^{2p}p!(h)^p+&&\\
&\quad +8p2^{2p-2}(p-1)![h]^{p-1}(h+p-1) + O(x^2)&&\text{by \eqref{y2p}}\\
&=(1-{1\over 2}[x,y]+{1\over 2}yx)2^{2p+1}p![h]^p+O(x^2)&&\\
&=(1-{1\over 2}h+{1\over 2}yx)2^{2p+1}p![h]^p+O(x^2)&&\text{by \eqref{xy}}\\
&=-2^{2p+1}p![h-1]^{p+1}+2^{2p}p![h]^pyx+O(x^2)&&\text{by \eqref{ypqh}, \eqref{xpqh}}
\end{align*}
as required.
\end{proof}

\ssection{Structure of the spin bundle}\label{spinstructure}

\begin{proposition}\label{extension}  Let $\psi\in\Gamma(\SSS(\tilde{G}))$.  The following are equivalent
\begin{enumerate}
\item $\psi$ extends some spinor $\Psi\in\Gamma(\SSS(G))$.
\item $(x\psi)|_G=0$.
\item There is a $\phi\in\Gamma(\SSS(\tilde{G}))$ such that $\psi=x\phi$.
\end{enumerate}
\end{proposition}

First, there is a useful lemma, which is also of independent interest:

\begin{lemma}\label{ker=im} Regard $x$ as a linear mapping $\SSS(\tilde{G})|_G\rightarrow \SSS(\tilde{G})|_G$.  Then the kernel of $x$ is equal to the image of $x$.  
\end{lemma}
\begin{proof}
By \eqref{x|_G}, $x^2=0$.  Hence $\im x\subseteq \ker x$.  For the reverse inclusion, it suffices to prove that $\dim\ker x=\dim\im x$.

Note that the Euler field $x^a$ is null on $G$ by \eqref{x|_G}.  Therefore, after a dilation as necessary, the vector $x^a$ can be represented with respect to an orthonormal basis by
$$x^0=1, x^1=1, x^i=0\hbox{\ \ for\ }i\ge 2.$$
Thus, $x=\gamma_0+\gamma_1$ where $(\gamma_0)^2=1$, $(\gamma_1)^2=-1$.  The condition $x\psi=0$ on $G$ is equivalent to
\begin{equation}\label{kill1}(\gamma_0+\gamma_1)\psi=0.\end{equation}
Multiplying \eqref{kill1} by $\gamma_0$ gives
\begin{equation}\label{kill2}(1+\gamma_0\gamma_1)\psi=0.\end{equation}
Let $T=\gamma_0\gamma_1$.  Then
\begin{equation*}
T^2=\gamma_0\gamma_1\gamma_0\gamma_1=-(\gamma_0)^2(\gamma_1)^2=1.
\end{equation*}
Thus the only eigenvalues of $T$ are $\pm 1$.  It suffices to prove that $T$ has an equal number of $+1$ and $-1$ eigenvalues, for this then implies that $\dim\ker x=\dim\im x$.
  
We have that
\begin{align*}
-T&=\gamma_1\gamma_0\\
&=-\gamma_1\gamma_0\gamma_1\gamma_1\\
&={\gamma_1}^{-1}T\gamma_1.
\end{align*}
Hence $T$ is similar to its negative, and therefore has the same number of $+1$ and $-1$ eigenvalues.
\end{proof}

\begin{proof}[Proof of Proposition.]  

$1\Longrightarrow 2$:  Let $\Psi$ be a spinor in $\Gamma(\SSS(G))$, and let $\psi$ be an extension.  Then $(x\psi)|_G=x|_G \Psi=0$ by Observation \ref{SGstruct1}.  

$2\Longrightarrow 1$:    It follows from \eqref{SGstruct2} and the general theory of spinors that the fibre dimension of $\SSS(\tilde{G})$ is twice the fibre dimension of $\SSS(G)$.  Hence, by Lemma \ref{ker=im}, $\dim\ker x|_G=\dim \SSS(G)_p$.  But, because of the proof of $1\Rightarrow 2$, $\SSS(G)\subseteq \ker x$; hence we must have $\ker x=\SSS(G)$.

$2\Longleftrightarrow 3$:  Follows immediately from Lemma \ref{ker=im}.
\end{proof}

\begin{definition}
A spinor $\psi\in\Gamma(\SSS(\tilde{G}))$ is $O(x^k)$ if there exists $\theta\in\Gamma(\SSS(\tilde{G}))$ such that
$$\psi=x^k\theta.$$
\end{definition}

The fact that we have used $\tilde{G}$ here makes a subtle but important difference.  (For instance, every section of $\SSS(G)$ is $O(x)$ according to our definition, by Proposition \ref{extension}.)  Along these lines, it is convenient to introduce some more notation.

\begin{definition}  Let $p$ be any integer, let $q$ be nonnegative integer with $q>p$, and let $w\in\R$.  Define
\begin{align*}
\sss(p,\infty)[w]&=x^p\SSS(\tilde{G})[w-p]&&\text{for $p\ge 0$,}\\
\sss(p,\infty)[w]&=\SSS(\tilde{G})[w]&&\text{for $p<0$}
\end{align*}
and let
$$\sss(p,q)[w]=\sss(p,\infty)[w]/s(q,\infty)[w].$$
\end{definition}

Thus, $\sss(p,q)$ is spinors of the form $x^p\psi$ modulo $O(x^q)$.  Note that $\SSS(G)[w]$ is identified with $\sss(1,2)[w]$, and that $\SSS(\tilde{G})[w]|_G$ is identified with $\sss(0,2)[w]$.  Every arrow in the sequence 
$$\sss(0,q)[w]\xrightarrow{x}\sss(1,q+1)[w+1]\xrightarrow{x}\sss(2,q+2)[w+2]\xrightarrow{x}\dots$$
is an isomorphism, although it is natural to maintain a distinction between the individual members of the above sequence, as they are different components of an algebraic system.

\begin{lemma}\label{xinverse}  For $p,q$ any nonnegative integers with $q>p$,
\begin{equation}\label{whatydoes}y:\sss(p,q)[w]\rightarrow\sss(p-1,q-1)[w-1].\end{equation}
Furthermore, if $p>0$, $q>p+1$, and 
$$w\not\in\biggl\{\lfloor {p\over 2}\rfloor-{n\over 2},\lfloor {p\over 2}\rfloor-{n\over 2}+1,\dots,\lfloor{q\over 2}\rfloor-{n\over 2}-1\biggr\}$$
then $y$ is an isomorphism.  In the exceptional case $p>0$ and $q=p+1$, there are two subcases:

(a) $p=2k$ is even.  In this case, $y$ is always an isomorphism, and furthermore
\begin{equation}\label{xinverseeven}
-{1\over p}y=x^{-1}:\sss(p,p+1)[w]\rightarrow\sss(p-1,p)[w-1]
\end{equation}

(b) $p=2k+1$ odd.  In this case, $y$ is an isomorphism provided that $w\not=k-{1\over 2}n-1$, and
\begin{equation}\label{xinverseeq}
{1\over 2}(h-k)^{-1}y=x^{-1}:\sss(p,p+1)[w]\rightarrow\sss(p-1,p)[w-1]
\end{equation}
\end{lemma}
\begin{proof}  
Consider for instance
\begin{align*}
y(x^{2k}\psi_1+x^{2k+1}\psi_2)&=-2kx^{2k-1}\psi_1+x^{2k}y\psi_1+&&\\
&\quad +2x^{2k}(h+k)\psi_2-x^{2k+1}\psi_2&&\text{by \eqref{x2py} and \eqref{x2p+1y}},
\end{align*}
and \eqref{whatydoes} follows by applying such a calculation successively.  To show that this is an isomorphism, suppose given an element of $\sss(2k-1,q-1)[w-1]$:
\begin{equation*}
\phi=x^{2k-1}\phi_1+x^{2k}\phi_2+\dots
\end{equation*}
Given that $x^{2k-1}\phi_1=2kx^{2k-1}\psi_1$ one knows $\psi_1$.  If $(h+k)\psi_2\not=0$, then given $\psi_1$ and $\phi_2$, one solves the following equation
for $\psi_2$:
\begin{equation*}
x^{2k-2}\phi_2=y\psi_1+2x^{2k}(h+k)\psi_2.
\end{equation*}
In the general case, one can always group terms according to their powers of $x$, and inductively solve a system of linear equations, assuming that everything is of the appropriate homogeneity.

Now for the exceptional case $p>0$ and $q=p+1$, suppose given 
$$x^p\psi_1\in\sss(p,p+1)[w].$$  
If $p=2k$, then
\begin{align*}
yx^{2k}\psi&=-2kx^{2k-1}\psi+x^{2k}y\psi&&\text{by \eqref{x2py}}\\
&=-2kx^{2k-1}\psi+O(x^{2k})&&
\end{align*}
which proves \eqref{xinverseeven}.  If $p=2k+1$, then
\begin{align*}
yx^{2k+1}\psi&=2(h-k)x^{2k}\psi+x^{2k}y\psi&&\text{by \eqref{x2p+1y} and \eqref{xpqh}}\\
&=2(h-k)x^{2k}\psi+O(x^{2k+1})&&
\end{align*}
which proves \eqref{xinverse}.
\end{proof}

\ssection{Invariant powers of the Dirac operator: Even case}\label{evencase}

\begin{theorem}  Let $\Psi\in\sss(1,2)[w]$.

(a)  If $w+{1\over 2}n$ is not a positive integer, then there exists a unique representative $\psi$ of $\Psi$ in $\sss(1,\infty)[w]$ such that $y^2\psi=0$.

(b)  If $w=p-{1\over 2}n$ with $p$ a positive integer, then there exists a unique representative $\psi$ of $\Psi$ in $\sss(1,2p)[w]$ such that $y^2\psi=0$.  The obstruction at order $2p$ is given by $(x^{-2p+2}y^2\psi)\mod O(x^2)$ and is independent of the choice of representative $\psi$.
\end{theorem}
\begin{proof}
Let $\psi_0$ be an arbitrary representative of $\Psi$ in $\sss(1,\infty)$.  Inductively set
$$\psi_p=\psi_{p-1}+x^{2p}\phi$$
for some $\phi\in\sss(0,\infty)[w-2p]$, and solve the equation $y^2\psi_p=O(x^{2p})$.  One has
\begin{align*}
y^2\psi_p&=y^2\psi_{p-1}+y^2x^{2p}\phi&&\\
&=y^2\psi_{p-1}+4px^{2p-2}(h+p-1)\phi+x^{2p}y^2\phi&&\text{by \eqref{x2py2}}\\
&=y^2\psi_{p-1}+4px^{2p-2}(w+{1\over 2}n-p)\phi+O(x^{2p}).&&
\end{align*}
Thus the equation $y^2\psi_p=O(x^{2p})$ can be solved for $\phi$ uniquely modulo $O(x^2)$ provided that $w\not=p-{1\over 2}n$.  In case $w=p-{1\over 2}n$, $y^2\psi_p=y^2\psi_{p-1}$ modulo $O(x^{2p})$, so that $x^{-2p+2}y^2\psi_p\mod O(x^2)$ is independent of the choice of $\phi$, and lies in $\sss(0,2)$.
\end{proof}

More important for the present purposes is not the obstruction $x^{-2p+2}y^2\psi \mod O(x^2)$, but rather $x^{-2p+3}y^2\psi \mod O(x^2)$, since the latter is evidently in $\sss(1,2)$.  Furthermore $x^{-2p+3}y^2\psi \mod O(x^2)$ gives the obstruction at the top degree in $x$.   The mapping
$$L_{2p}:\Psi\mapsto x^{-2p+3}y^2\psi \mod O(x^2)$$
defines a conformally invariant operator $\sss(1,2)[p-{1\over 2}n]\rightarrow \sss(1,2)[-p-{1\over 2}n+1]$.

\begin{theorem}  Let $p$ be an integer, let $w=p-{1\over 2}n$, let $\Psi\in\sss(1,2)[w]$, and let $\psi\in\sss(1,\infty)[w]$ be any representative of $\Psi$.  Then $y^{2p}\psi\mod O(x^2)\in\sss(0,2)[w]$ does not depend on the choice of representative.  Thus 
$$R_{2p}:\Psi\mapsto xy^{2p}\psi\mod O(x^2)$$
is a conformally invariant operator $\sss(1,2)[p-{1\over 2}n]\rightarrow\sss(1,2)[-p-{1\over 2}n+1]$.
\end{theorem}

\begin{proof}  It suffices by Proposition \ref{extension} to check that $y^{2p}x^2\theta =0 \mod O(x^2)$ for all $\theta\in\sss(0,\infty)[w-2]$.  Thus consider
\begin{align*}
y^{2p}x^2\theta&=[y^{2p},x^2]\theta+O(x^2)&&\\
&=4py^{2p-2}(h-p+1)\theta+O(x^2)&&\\
&=O(x^2)&&\text{since $\theta\in\sss(0,\infty)[w-2]$}
\end{align*}
\end{proof}

\begin{theorem}  $L_{2p}$ is a nonzero multiple of $R_{2p}$.
\end{theorem}
\begin{proof}
Let $w=p-{1\over 2}n$, and $\Psi\in\sss(1,2)[w]$.  Let $\psi$ be the representative of $\Psi$ with $y^2\psi=O(x^{2p-2})$.  Then $y^2\psi=x^{2p-3}L_{2p}\Psi$.  Now multiply both sides by $y^{2p-2}$ to give
\begin{equation}\label{y2ppsi}
y^{2p}\psi=y^{2p-2}x^{2p-3}L_{2p}\Psi.
\end{equation}
Write $L_{2p}\Psi=x\beta$ for some $\beta\in\sss(0,1)[w-2p]$.  Then \eqref{y2ppsi} gives
\begin{align*}
y^{2p}\psi&=y^{2p-2}x^{2p-2}\beta&&\\
&=2^{2p-2}(p-1)![h]^{p-1}\beta+O(x^2)&&\text{by \eqref{y2p}}\\
&=(-1)^{p-1}2^{2p-2}(p-1)!^2\beta+O(x^2)&&\text{by homogeneity}
\end{align*}
Thus
\begin{align*}
xy^{2p}&=(-1)^p2^{2p-2}(p-1)!^2x\beta+O(x^2)\\
&=(-1)^p2^{2p-2}(p-1)!^2L_{2p}\Psi+O(x^2),
\end{align*}
which is the required result.
\end{proof}

\ssection{Invariant powers of the Dirac operator: Odd case}\label{oddcase}

\begin{theorem}\label{confdirac2}
Let $\Psi\in\sss(1,2)[w]$.

(a)  If $w+{1\over 2}n-1$ is not a positive integer, then there exists a unique representative $\psi$ of $x^{-1}\Psi$ in $\sss(0,\infty)[w-1]$ such that $y\psi=0$.

(b)  If $w=p-{1\over 2}n+1$ with $p$ a positive integer, then $x^{-1}\Psi$ has a unique representative $\psi\in\sss(0,2p+1)[w-1]$ which obeys $y\psi=0$.  The obstruction to extension of $\psi$ at order $O(x^{2p+1})$ is given by $(x^{-2p+1}y\psi)\mod O(x^2)$, which is independent of the choice of representative $\psi$ and therefore defines a conformally invariant operator on $\sss(1,2)[w]$.
\end{theorem}

\begin{proof}  Let $\psi_0\in\sss(0,\infty)$ be an arbitrary representative of $x^{-1}\Psi$.  Inductively, write
$$\psi_k=\psi_{k-1}+x^k\phi$$
for some spinor $\phi\in\sss(0,\infty)[w-k-1]$, and suppose that $y\psi_{k-1}=O(x^{k-1})$.  One may now solve for $\phi$ such that $y\psi_k=O(x^k)$.  There are two distinct cases:

{\it Case 1.} $k=2p$.
\begin{align*}
y\psi_{2p}&=y\psi_{2p-1}+yx^{2p}\phi&&\\
&=y\psi_{2p-1}-[x^{2p},y]\phi+O(x^{2p})&&\\
&=y\psi_{2p-1}-2px^{2p-1}\phi+O(x^{2p})&&\text{by \eqref{x2py}}\\
&=x^{2p-1}\beta-2px^{2p-1}\phi+O(x^{2p})&&\text{for some $\beta$}
\end{align*}
(where the last step uses the induction hypothesis).
So setting $\phi={1\over 2p}\beta$ gives the desired result.

{\it Case 2.} $k=2p+1$.
\begin{align*}
y\psi_{2p+1}&=y\psi_{2p}+yx^{2p+1}\phi&&\\
&=y\psi_{2p}-[x^{2p+1},y]\phi+O(x^{2p+1})&&\\
&=y\psi_{2p}-2x^{2p}(h+p)\phi+O(x^{2p+1})&&\text{by \eqref{x2p+1y}}\\
&=x^{2p}\beta-2x^{2p}(w-1+{n\over 2}-p)\phi+O(x^{2p+1})&&\text{by induction.}
\end{align*}
The equation $y\psi_{2p+1}=O(x^{2p+1})$ can be solved for $\phi$ if and only if $w\not=p-{n\over 2}+1$.  In case $w=p-{n\over 2}+1$, $y\psi_{2p}=y\psi_{2p+1}$ modulo $O(x^{2p+1})$ independently of $\phi$.  Hence, $x^{-2p+1}y\psi_{2p+1}$ is independent of $\phi$ modulo $O(x^2)$.  
\end{proof}

Thus, to summarize, in the notation of the preceding theorem, the mapping
$$\Psi\mapsto L_{2p+1}\Psi\overset{def}{=}(x^{-2p+1}y\psi_{2p+1})\mod O(x^2)$$
is a conformally invariant operator.  Furthermore, $L_{2p+1}$ maps into $\sss(1,2)$ since
\begin{align*}
xL_{2p+1}\Psi\mod O(x^2)&=x^{-2p+2}y\psi_{2p+1}\mod O(x^2)&&\\
&=x^{-2p+2}y\psi_{2p}\mod O(x^2)&&\\
&=(x^{-2p+2}O(x^{2p}))\mod O(x^2)&&\\
&=0
\end{align*}
Therefore,
$$L_{2p+1}:\sss(1,2)[p-{n\over 2}+1]\rightarrow\sss(1,2)[-p-{n\over 2}].$$

For the next theorem, given $\Psi\in\sss(1,2)[w]$, one needs to be able to define a distinguished class of representatives in $\sss(1,3)[w]$.  This distinguished class is supplied by the following.  Let $\Psi\in\sss(1,2)[w]$.  Let $\psi_0\in\sss(0,2)[w-1]$ be the unique representative of $x^{-1}\Psi$ such that $y\psi_0=0$.  Then $x\psi_0\in\sss(1,3)[w]$ is a uniquely determined representative of $\Psi$.  Call this representative the {\it preferred} representative.

\begin{theorem}  Let $\Psi\in\sss(1,2)[p-{1\over 2}n+1]$, let $\bar{\psi}$ be the preferred representative in $\sss(1,3)$, and let $\psi$ be any representative of $\bar{\psi}$ in $\sss(1,\infty)$.  Then $y^{2p+1}\psi\mod O(x^2)$ lies in $\sss(1,2)[-p-{1\over 2}n]$ and is independent of the $\psi$ chosen to represent $\bar{\psi}$.  Therefore it defines a conformally invariant operator
$$R_{2p+1}:\sss(1,2)[p-{1\over 2}n+1]\rightarrow\sss(1,2)[-p-{1\over 2} n].$$
\end{theorem}

\begin{proof}  To see that $y^{2p+1}\psi\mod O(x^2)$ does not depend on the extension $\psi$, it suffices to check that $y^{2k+1}x^3\theta\mod O(x^2)=0$ for all $\theta\in\sss(0,\infty)[p-{1\over 2}n-2]$. 
Hence consider
\begin{align*}
y^{2p+1}x^3\theta&=[y^{2p+1},x^2]x\theta+x^2y^{2p+1}x\theta&&\\
&=4py^{2p-1}(h-p)x\theta+2y^{2p}x^2\theta+x^2y^{2p+1}x\theta&&\text{by \eqref{y2p+1x2}}\\
&=4py^{2p-1}x(h-p+1)\theta+2[y^{2p},x^2]\theta+O(x^2)&&\\
&=4py^{2p-1}x(h-p+1)\theta+8py^{2p-2}(h-p+1)\theta+O(x^2)&&\text{by \eqref{y2px2}}\\
&=4p(y^{2p-1}x+2y^{2p-2})(h-p+1)\theta+O(x^2)&&\\
&=O(x^2)&&\text{by homogeneity}
\end{align*}
as required.  So $R_{2p+1}$ is a well-defined operator, and it maps into $\sss(0,2)[-p-{1\over 2}n]$.  

Now we must check that $R_{2p+1}\Psi\in\sss(1,2)$, i.e., that $xR_{2p+1}\Psi=0$.
\begin{align*}
xy^{2p+1}x\phi&=x[y^{2p+1},x]\phi+O(x^2)&&\\
&=2pxy^{2p}(h-p)\phi+O(x^2)&&\text{by \eqref{y2p+1x}}\\
&=O(x^2)&&\text{by homogeneity.}
\end{align*}
\end{proof}

\begin{theorem}\label{confdirac3}  $L_{2p+1}$ is a nonzero multiple of $R_{2p+1}$.
\end{theorem}

\begin{proof}  Let $\psi$ be a representative of $x^{-1}\Psi$ in $\sss(0,2p+1)[w-1]$ such that $y\psi=0$.  Then $y\psi=x^{2p-1}L_{2p+1}\Psi+O(x^{2p+1})$.  Note also that 
\begin{equation}\label{psiphi}\psi\mod O(x^2)={1\over 2}h^{-1}y\phi\end{equation}
for some $\phi\in\sss(1,3)[w]$ by \eqref{xinverseeq}, $\phi$ being a preferred representative of $\Psi$.  On the other hand
$${1\over 2}h^{-1}y\psi=Cy\psi$$
for some nonzero constant $C$, since $\psi$ is homogeneous.  Multiplying both sides of \eqref{psiphi} by $y^{2p-1}$ gives
\begin{align*}
C y^{2p+1}\phi&=y^{2p}\psi&&\\
&=y^{2p-1}x^{2p-1}L_{2p+1}\Psi+O(x^2)&&\\
&=-2^{2p-1}(p-1)![h-1]^pL_{2p+1}\Psi+&&\\
&\quad +2^{2p-2}(p-1)![h]^{p-1}yxL_{2p+1}\Psi+O(x^2)&&\text{by \eqref{y2p+1}}\\
&\overset{def}{=}(\text{term 1})+(\text{term 2})+O(x^2)&&
\end{align*}
Throughout the remainder of this calculation, terms which are $O(x^2)$ will be dropped without comment.  For (term 1), note that $hL_{2p+1}\Psi=-(p-1)$, and so $[h-1]^pL_{2p+1}\Psi=(-1)^pp!\Psi$.  Thus
\begin{equation}\label{term1}
(\text{term 1})=(-1)^{p-1}2^{2p-2}(p-1)!p!L_{2p+1}\Psi
\end{equation}

For (term 2), let $A=2^{2p-2}(p-1)![h]^{p-1}$ and $B=2AhL_{2p+1}\Psi$.  Then 
\begin{align*}
(\text{term 2})&=AyxL_{2p+1}\Psi&&\\
&=B-AxyL_{2p+1}\Psi&&\\
&=B-Axyx^{-2p+1}y\psi&&\\
&=B-Ax^{-2p+2}y^2\psi+2Ax^{-2p+1}(h-p)y\psi&&\text{by \eqref{x2p+1y}}\\
&=B-Ax^{-2p+2}y^2\psi+2Ax^{-2p+1}y(h-p-1)\psi&&\\
&=B-Ax^{-2p+2}y^2\psi&&\text{by homogeneity}
\end{align*}

The factor $x^{-2p+2}y^2\psi$ can be rewritten as follows.  Let $\beta$ be a spinor such that $x\beta=L_{2p+1}\Psi$.  Then
\begin{align*}
x^{-2p+2}y^2\psi&=x^{-2p+2}yx^{2p}\beta&&\\
&=x^2y\beta-2px\beta&&\text{by \eqref{x2py}}\\
&=-2pL_{2p+1}\Psi&&
\end{align*}
Thus
\begin{align*}
\text{(term 2)}&=B+2pAL_{2p+1}\Psi&&\\
&=(2^{2p-1}(p-1)!h[h]^{p-1}+2^{2p-1}p![h]^{p-1})L_{2p+1}\Psi&&\\
&=2^{2p-1}\bigl((-1)^{p-1}(p-1)!+(-1)^pp!+(-1)^{p-1}p!\bigr)p!L_{2p+1}\Psi&&\\
&=(-1)^{p-1}2^{2p-1}(p-1)!p!L_{2p+1}\Psi&&
\end{align*}
Hence $\text{(term 1)} + \text{(term 2)}$ is a nonzero multiple of $L_{2p+1}\Psi$, as required.
\end{proof}

\end{document}